\documentclass[12pt]{article}
\usepackage{amsmath, amsthm, amsfonts}

\theoremstyle{plain}
\newtheorem{theorem}{Theorem}

\newtheorem{prop}[theorem]{Proposition}

\newcommand{\p}{\ensuremath{\mathbf{p}}}
\newcommand{\q}{\ensuremath{\mathbf{q}}}
\newcommand{\ve}{\ensuremath{\mathbf{v}}}
\newcommand{\norm}[1]{\ensuremath{\left\Vert #1 \right\Vert}}

\newcommand{\reals}{\ensuremath{\mathbb{R}}}
\newcommand{\natnum}{\ensuremath{\mathbb{N}}}
\newcommand{\infnorm}[1]{\ensuremath{\left\Vert #1 \right\Vert_\infty}}

\DeclareMathOperator{\mat}{Mat}
\DeclareMathOperator{\GCD}{GCD}

\begin{document}

\title{A quantitative Khintchine--Groshev type theorem over a field of
  formal series}
\author{M. M. Dodson$^1$, S. Kristensen$^2$ and J. Levesley$^3$\\[2pt] 
{\small $^1$Department of Mathematics, University of York,} \\
{\small Heslington, York, YO10 5DD, UK}\\ 
{\small \texttt{mmd1@york.ac.uk}}\\[2pt]
{\small $^2$School of Mathematics, University of Edinburgh, JCMB,}\\  
{\small King's Buildings, Mayfield Road, Edinburgh, EH9 3JZ, UK} \\
{\small \texttt{Simon.Kristensen@ed.ac.uk}}\\[2pt]
{\small $^3$Department of Mathematics, University of York,} \\
{\small Heslington, York, YO10 5DD, UK}\\ 
{\small \texttt{jl107@york.ac.uk}}}

\maketitle
\thispagestyle{empty}

\begin{abstract}
  An asymptotic formula which holds almost everywhere is obtained for
  the number of solutions to the Diophantine inequalities $\norm{\q
    A-\p}<\psi(\norm{\q})$, where $A$ is an $n \times m$ matrix ($m >
  1$) over the field of formal Laurent series with coefficients from a
  finite field, and $\p$ and $\q$ are vectors of polynomials over the
  same finite field. 
\end{abstract}

\section{Introduction}
\label{sec:introduction}

Let $\mathbb{F}$ denote the finite field of $k = p\,^l$ elements, where
$p$ is a prime and $l$ is a positive integer. We define
\begin{equation}
  \label{eq:1}
  \mathcal{L} = \left\{\sum_{i=-n}^\infty a_{-i} X^{-i} : n \in
    \mathbb{Z}, a_i \in \mathbb{F}, a_n \neq 0\right\} \cup \{0\}.
\end{equation}
Under usual addition and multiplication, this set is a field,
sometimes called \emph{the field of formal Laurent series with
  coefficients from $\mathbb{F}$}. We may define an absolute value on
$\mathcal{L}$ by setting
\begin{displaymath}
  \norm{\sum_{i=-n}^\infty a_{-i} X^{-i}} = k^n, \quad \norm{0} = 0.
\end{displaymath}
This absolute value is ultra-metric. Under the induced metric,
$d(x,y) = \norm{x-y}$, the space $(\mathcal{L}, d)$ is a complete
metric space.

The approximation of elements of $\mathcal{L}$ by ratios of elements
in the polynomial ring $\mathbb{F}[X]$ has been studied extensively
(see \emph{e.g.}~the survey papers by Lasjaunias \cite{MR2001k:11135}
and Schmidt \cite{MR2001j:11063}) and has been used in the analysis of
pseudorandom sequences employed in cryptography by Niederreiter and
Vielhaber~\cite{nv1997}.

In this paper, we are concerned with the metrical theory of such
Diophantine approximations. Let $\psi: \reals_+ \rightarrow \reals_+$
be a function with $\psi(x)$ non-increasing.  In \cite{MR43:161}, de
Mathan showed that the set of elements $x \in \mathcal{L}$ for which
the inequality
\begin{displaymath}
  \norm{qx-p} < \psi (\norm{q})
\end{displaymath}
has infinitely many solutions $q,p \in \mathbb{F}[X]$, $q \neq 0$ is
null or full (with respect to the Haar measure) accordingly as the
series $\sum_{r=1}^\infty \psi(r)$ diverges or converges. This was
extended to systems of linear forms in Kristensen \cite{kristensen03},
as follows.
\begin{theorem}[\protect{\cite[Theorem 3]{kristensen03}}]
  Let $\psi: \reals_+ \rightarrow \reals_+$ be decreasing. Let $m, n
  \in \natnum, m \geq 2$. The set of $m \times n$ matrices $A$ with
  entries from $\mathcal{L}$ for which the inequalities
  \begin{equation}
    \label{eq:2}
    \infnorm{\mathbf{q}A-\mathbf{p}} < \psi(\infnorm{\mathbf{q}})
  \end{equation}
  have infinitely many solutions $\mathbf{p} \in \mathbb{F}[X]^n$,
  $\mathbf{q} \in \mathbb{F}[X]^m, \mathbf{q} \neq \mathbf{0}$ is null
  or full accordingly as the series $\sum_{r=1}^\infty r^{m-1}
  \psi(r)^n$ converges or diverges, where $\infnorm{\mathbf{q}} = \max
  \{ \norm{q_i}\}$ for $\mathbf{q} = (q_1, \dots, q_m)$.
\end{theorem}

Here, we are concerned with the asymptotic number of solutions to the
inequalities \eqref{eq:2}. In the real case, the analogous asymptotics
were found by Schmidt, first in the case of simultaneous approximation
as well as approximation of a single linear form in Schmidt
\cite{MR22:9482} and since for systems of linear forms as well as for
restricted sets of $\mathbf{q}$'s  in Schmidt \cite{MR28:3018}.

We will restrict ourselves to considering error functions taking their
values in the set $V = \{k^{-n} : n \in \mathbb{N}\}$. For general
error functions, see the remark following the statement of the
theorem. We will prove the following theorem:
\begin{theorem}
  \label{thm:main_theorem}
  Let $\epsilon > 0$, let $\psi : \mathbb{R}_+ \rightarrow V$ and let
  $N(Q,A)$ denote the number of solutions to \eqref{eq:2} with
  $\infnorm{\mathbf{q}} \leq k^Q$. Let
  \begin{displaymath}
    \Phi(Q) =
      m (k-1) k^{m-1} \sum_{r=0}^Q k^{rm} \psi(k^r)^n
  \end{displaymath}
  Then
  \begin{displaymath}
    N(Q,A) = \Phi(Q) + O\left(\Phi(Q)^{1/2}\log^{3/2 +
        \epsilon}\left(\Phi(Q)\right)\right)
  \end{displaymath}
  for almost every $m\times n$ matrix $A$ with entries from
  $\mathcal{L}$.
\end{theorem}

The reason for restricting the choice of error functions is that the
only possible distances in the space $\mathcal{L}$ are of the form
$k^r$ where $r \in \mathbb{Z}$. For other error functions, we could
define a function, $\lfloor \cdot \rfloor : \mathbb{R}_+ \rightarrow
V$ say, mapping $x \in \mathbb{R}$ to the unique number $\lfloor x
\rfloor \in V$ such that $\lfloor x \rfloor \leq x < k \lfloor x
\rfloor$. On replacing $\psi(\cdot)$ with $\lfloor \psi(\cdot)
\rfloor$ at every occurence, we would obtain the theorem for general
decreasing error functions.  However, for ease of notation we consider
only the restricted case.

\section{Proof of main theorem}
\label{sec:proof-main-theorem}

The proof has two main ingredients. The first has to do with the
geometry of the underlying vector spaces. The second is a purely
probabilistic theorem. We first prove the geometrical results.

We identify $\mat_{m\times n}(\mathcal{L})$ with $\mathcal{L}^{mn}$.
Define for any $\q \in \mathbb{F}[X]^m$ the set
\begin{equation}
  \label{eq:6}
  B_\q = \left\{A \in I^{mn} : \inf_{\p \in \mathbb{F}[X]^n}
  \infnorm{\q A - \p} < \psi(\infnorm{\q})\right\},
\end{equation}
where $I^{mn}$ denotes the $\infnorm{\cdot}$-unit ball in
$\mathcal{L}^{mn}$. The Haar measure on $\mathcal{L}^{mn}$, normalised
so that the measure of $I^{mn}$ is equal to $1$, will be denoted by
$\mu$.

We will prove the following propositions:
\begin{prop}
  \label{prop:invariance}
  \begin{displaymath}
    \mu\left(B_\q)\right) = \psi(\infnorm{\q})^n. 
  \end{displaymath}
\end{prop}

\begin{prop}
  \label{prop:independence}
  Let $\q, \q' \in \mathbb{F}[X]^m$ be linearly independent over
  $\mathcal{L}$. Then
  \begin{displaymath}
    \mu\left(B_\q \cap B_{\q'} \right) = \mu\left(B_\q \right)
    \mu\left(B_{\q'}\right). 
  \end{displaymath}
\end{prop}

In both proofs, we follow the method from Dodson \cite{MR95d:11092}. 

\begin{proof}[Proof of Proposition \ref{prop:invariance}]
  By the rank equation, the solution curves to the equations $\q A =
  \p$ are $(m-1)n$ dimensional affine spaces over $\mathcal{L}$. We
  begin by calculating the number of affine spaces which pass through
  the unit ball. First, note that if there is a solution to the
  equation $\q A = \p$ with $A \in I^{mn}$, then
  \begin{equation}
    \label{eq:7}
    \infnorm{\p} = \infnorm{\q A} \leq \infnorm{\q}\infnorm{A} <
    \infnorm{\q},
  \end{equation}
  so certainly, the condition $\infnorm{\p} < \infnorm{\q}$ is
  necessary. We claim that it is also sufficient.
  
  For this, it suffices to find a solution $A \in I^{mn}$ which
  satisfies the equation. Suppose that $\infnorm{\p} < \infnorm{\q}$.
  We assume without loss of generality that $\infnorm{\q} =
  \norm{q_1}$.  Now,
  \begin{equation}
    \label{eq:8}
    \q A  = \q
    \begin{pmatrix}
      p_1/q_1 & \cdots & p_n/q_1 \\
      0 & \cdots & 0 \\
      \vdots & & \vdots \\
      0 & \cdots & 0
    \end{pmatrix} = \p
  \end{equation}
  and $A \in I^{mn}$.
  
  As in Dodson \cite{MR95d:11092}, we consider the simplest
  non-trivial case where $\q=(q_1,q_2)$ and $\p=p$ and subsequently
  extend this to the general case. In this case, the solution curves
  to the equations $\q A = p$ define $\infnorm{\q}$ affine
  $1$-dimensional spaces in $I^2$.  These partition $I^2$ into
  $\infnorm{\q}$ strips, $\tilde{S}_i$ say, defined by inequalities
  $\norm{\q A - p} < 1$. The measure of each such strip may be
  calculated using a characterisation of a translation invariant
  measure due to Mahler (see \cite{MR2:350c}), which implies
  that the measure of a parallelogram is $1 / \det(w_1, w_2)$, where
  $w_1$ and $w_2$ are the spanning vectors. Since the distance between
  each affine $1$-space is $1/\infnorm{\q}$, the solution curves
  partition $I^2$ into sets of the same size, $\mu (\tilde{S}_i) =
  1/\infnorm{\q}$. By the same characterization, we find that around
  each solution curve we have a component, $B_i$ say, of the set $B_\q$
  of measure $\psi(\q)/\infnorm{\q}$.  Hence
  \begin{equation}
    \label{eq:9}
    \mu\left(B_\q\right) =
    \dfrac{\mu\left(B_\q\right)}{\mu(I^2)} 
    = \dfrac{\mu\left(\cup B_i\right)}{\mu\left(\cup
        \tilde{S}_i\right)} = \dfrac{\mu(B_i)}{\mu(\tilde{S}_i)} = 
    \dfrac{\psi(\infnorm{\q})/\infnorm{\q}}{1/\infnorm{\q}} =
    \psi(\infnorm{\q}). 
  \end{equation}

  To obtain the proposition for general $m,n \in \natnum$, consider
  $n$ copies of the span of $\q$ and apply the above argument to
  resulting prisms in $I^{mn}$. This implies the proposition.
\end{proof}

\begin{proof}[Proof of Proposition \ref{prop:independence}]
  Again, we consider the simplest non-trivial case, $m=2, n=1$. Let
  $\q, \q' \in \mathbb{F}[X]^2$ be linearly independent. We calculate
  the number of intersections between the solution curves to the
  equations $\q A = p$ and the equations $\q' A = p'$, where $p, p'$
  runs over the possible values. This amounts to solving the system
  \begin{displaymath}
    \begin{pmatrix}
      q_1 & q_2 \\
      q_1' & q_2'
    \end{pmatrix}
    \begin{pmatrix}
      a_1 \\
      a_2
    \end{pmatrix}
    =
    \begin{pmatrix}
      p \\
      p'
    \end{pmatrix},
    \quad \infnorm{p} < \infnorm{\q}, \infnorm{p'} < \infnorm{\q'}.
  \end{displaymath}
  There are exactly $\norm{\det\bigl(\begin{smallmatrix} q_1 & q_2 \\
      q_1' & q_2' \end{smallmatrix}\bigr)}$ such solutions. To each
  such solution, we may assign a parallelogram defined by the
  inequality
  \begin{displaymath}
    \max\{\norm{\q A - p}, \norm{\q' A - p} \} < 1
  \end{displaymath}
  The parallelogram is seen to be of measure
  $1/\norm{\det\bigl(\begin{smallmatrix} q_1 & q_2 \\ q_1' & q_2'
    \end{smallmatrix}\bigr)}$, and the parallelograms are mutually
  disjoint. 

  To show that these parallelograms partition $I^2$, it remains to be shown that each of the parallelograms defined above
  is a proper subset of $I^2$. But this is the case, since any
  parallelogram may be written as
  \begin{displaymath}
    \left\{x \in \mathcal{L}^2 : x = \hat{q}_1 t_1 + \hat{q}_2 t_2 + p,
    t_1, t_2 \in I\right\}
  \end{displaymath}
  for some $\hat{q}_1, \hat{q}_2 \in \mathcal{L}^2$. Clearly,
  \begin{displaymath}
    \left\{x \in \mathcal{L}^2 : x = \hat{q}_1 t_1 + \hat{q}_2 t_2 + p,
    t_1, t_2 \in I\right\} 
  \subseteq B\left(p, \max\left(\infnorm{\hat{q}_1},
      \infnorm{\hat{q}_2}\right)\right), 
  \end{displaymath}
  so by the ultrametric property, the parallelogram is either fully
  contained in $I^2$ or disjoint with $I^2$. Since the parallelograms
  bounded by the solution curves are disjoint, there can be no more
  than the required number.

  Furthermore, around each intersection point, there is another
  parallelogram of measure
  $\psi(\q)\psi(\q')/\norm{\det\bigl(\begin{smallmatrix} q_1 & q_2 \\
      q_1' & q_2' \end{smallmatrix}\bigr)}$, constituting a part of
  $B_\q \cap B_{\q'}$ whenever it is a subset of $I^2$.

  With the above tools, we may apply a proportionality argument
  analogous to \eqref{eq:9} to obtain the proposition in this case.
  For the general case, we consider $n$ copies of the span of $\q$ and
  $\q'$ and apply the above to the $mn$ dimensional prisms to obtain
  the proposition.
\end{proof}

We are now ready to prove the main theorem.

\begin{proof}[Proof of Theorem \ref{thm:main_theorem}]
  The final ingredient in the proof is Lemma 10 in Sprind\-\v
  zuk~\cite{MR80k:10048}.
  
  Let $f_{\q}(A)$ be the characteristic function of $B_\q$, $f_{\q} =
  \psi(\norm{\q})^n$ and let $\tau(\q) = \psi(\q)^n d(\q)$, where
  $d(\q)$ denotes the number of common divisors in $\mathbb{F}[X]$ of
  the coordinates of $\q$.  Clearly, by
  Proposition~\ref{prop:invariance} and
  Proposition~\ref{prop:independence} for $s < t$,
  \begin{displaymath}
    \int \left(\sum_{k^s < \infnorm{\q} \leq k^t} f_\q(A) -
      \sum_{k^s < \infnorm{\q} \leq k^t} f_\q\right)^2 dA \ll
    \sum_{k^s < \infnorm{\q} \leq k^t} \tau_\q,
  \end{displaymath}
  as we only get contributions from the diagonal and elements that
  corresponding to pairs of parallel $\q$'s. By Lemma
  10 in Sprind\v zuk \cite{MR80k:10048}, we then have for almost every
  $A$, 
  \begin{equation}
    \label{eq:3}
    N(Q,A) = \sum_{\infnorm{\q} \leq k^Q} f_\q(A) 
    = \sum_{\infnorm{\q} \leq k^Q} f_\q + O
    \left(T(Q)^{1/2}\log^{3/2+\epsilon}T(Q)\right),
  \end{equation}
  where $T(Q) = \sum_{\infnorm{\q} \leq k^Q} \tau(\q)$. We need to
  prove that the right hand side is dominated by the first term.

  We begin with this term. By formula (1.4) in Kristensen
  \cite{kristensen03}, 
  \begin{displaymath}
    \sum_{\infnorm{\q} \leq k^Q} f_\q = \sum_{r=0}^Q
    \sum_{\infnorm{\q} = k^r} \psi(k^r)^n 
    = m (k-1) k^{m-1} \sum_{r=0}^Q k^{rm} \psi(k^r)^n = \Phi(Q).
  \end{displaymath}
  Thus we need only worry about the error term $T(Q)$.  Clearly, it
  suffices to prove that
  \begin{displaymath}
    T(Q) = O(\Phi(Q)).
  \end{displaymath}
  We first observe that by denoting $\q=(q_1, \dots, q_m)$ and again
  applying~(1.4) from Kristensen~\cite{kristensen03},
  \begin{equation*}
    \#\{\q\in f[X]^m\colon \infnorm{\q}=k^r\}= m(k-1)k^{m-1+rm},
  \end{equation*}
  which gives the number of $\q$ of given height $k^r$, we immediately
  obtain
  \begin{eqnarray*}
    T(Q) &=& \sum_{\infnorm{\q} \leq
      k^Q} \psi(\infnorm{\q})^n\sum_{d\vert(q_1,\dots,q_m)}1 \\ 
    &\ll & \sum_{r=0}^Q \sum_{\alpha = 0}^r
    \sum_{\substack{\infnorm{\q} = k^r\\
        \GCD(q_1,\dots,q_m)=a=k^\alpha } }  \psi(\infnorm{\q})^n 
    \sum_{d \vert a}\ 1 \\
    &\ll& \sum_{r=0}^Q \psi(k^r)^n\sum_{\alpha=0}^r
    \sum_{\substack{\infnorm{\ve} = k^{r-\alpha}\\
        \GCD(v_1,\dots,v_m)=1 } }1\\
    &\ll& \sum_{r=0}^Q \psi(k^r)^n\sum_{\alpha=0}^r
    m(k-1)k^{m-1}k^{(r-\alpha) m}\\
    &\ll&  m(k-1)k^{m-1} \sum_{r=0}^Q \psi(k^r)^nk^{rm}
    \sum_{\alpha=0}^r k^{-\alpha m}\\
    &\ll&   m(k-1)k^{m-1} \sum_{r=0}^Q \psi(k^r)^n k^{rm} \ll \Phi(Q), 
  \end{eqnarray*}
  as required.
\end{proof}

\section{Acknowledgements}
\label{sec:acknowledgements}

We thank the referee for some helpful comments.  Research partially
funded by EPSRC grant no. GR/N02832/01 with additional support from
INTAS grant no. 001--429.  SK is a William Gordon Seggie Brown Fellow.

\end{document}